\documentclass[3p]{elsarticle}

\usepackage{lineno,hyperref}
\modulolinenumbers[5]

\journal{Journal of \LaTeX\ Templates}

\usepackage{amssymb}
\usepackage{listings}
\usepackage{color}
\usepackage{float}
\usepackage{subfig}
\usepackage{algorithm}
\usepackage{algorithmic}
\usepackage[fleqn]{amsmath}
\usepackage{graphicx}
\usepackage{amsthm}
\usepackage{tikz}
\usetikzlibrary{arrows.meta}

\usepackage{array,booktabs,ragged2e}
\newcolumntype{P}[1]{>{\RaggedRight\arraybackslash}p{#1}}

\usepackage{booktabs} 

\newcommand{\V}[1]{\mathbf{#1}}
\newcommand{\M}[1]{\mathbf{#1}}
\newcommand{\Vg}[1]{\boldsymbol{#1}}
\newcommand{\tr}[1] {#1^{\sf{T}} }

\newtheorem{remark}{Remark}

\begin{document}

\begin{frontmatter}

\title{Fast multiscale contrast independent preconditioners for linear elastic topology optimization problems}





\author[4,1]{Miguel Zambrano}
\author[4,1]{Sintya Serrano}
\author[2,3]{Boyan S. Lazarov\corref{cor1}}
\ead{lazarov2@llnl.gov}
\author[1]{Juan Galvis}
\ead{jcgalvisa@unal.edu.co}
\cortext[cor1]{Corresponding author}
\address[1]{Departamento de Matem\'aticas, Universidad Nacional de Colombia, Carrera 45 No. 26-85, Edificio Uriel Gutierr\'ez, Bogot\'a D.C., Colombia}
\address[2]{Lawrence Livermore National Laboratory, Livermore, CA, US}
\address[3]{School of Mechanical, Aerospace and Civil Engineering, the University of Manchester, Manchester, UK}
\address[4]{Departamento de Ciencias Exactas, Universidad de las Fuerzas Armadas ESPE, Sangolqu\'i, Ecuador}

\begin{abstract}
The goal of this work is to present a fast and viable approach for the numerical solution of the high-contrast state problems arising in topology optimization. The optimization process is iterative, and the gradients are obtained by an adjoint analysis, which requires the numerical solution of large high-contrast linear elastic problems with features spanning several length scales. The size of the discretized problems forces the utilization of iterative linear solvers with solution time dependant on the quality of the preconditioner. The lack of clear separation between the scales, as well as the high-contrast, imposes severe challenges on the standard preconditioning techniques. Thus, here we propose new methods for the high-contrast elasticity equation with performance independent of the high-contrast and the multi-scale structure of the elasticity problem. The solvers are based on two-levels domain decomposition techniques with a carefully constructed coarse level to deal with the high-contrast and multi-scale nature of the problem. The construction utilizes spectral equivalence between scalar diffusion and each displacement block of the elasticity problems and, in contrast to previous solutions proposed in the literature, is able to select the appropriate dimension of the coarse space automatically. The new methods inherit the advantages of domain decomposition techniques, such as easy parallelization and scalability. The presented numerical experiments demonstrate the excellent performance of the proposed methods.
\end{abstract}


\end{frontmatter}


\section{Introduction} \label{s:intro}

Topology optimization \cite{Bendsoe2003} is an iterative design process aiming to find close to optimal material distribution by minimizing an objective function and fulfilling a set of constraints. More precisely, the discrete optimization problem, considered here, can be written as
\begin{align}
\label{eq1001}
{\rm{minimize}}_{\rho}: &\,\,g_0\left(\Vg{\rho}, \V{u}\right)= \V{f}^{\sf{T}}\V{u} \\
s.t.: &\,\, r\left(\Vg{\rho},\V{u}\right)=0, \quad \V{u}\in \mathcal{U}_{\text{ad}}\\
      &\,\, g_1\left(\Vg{\rho}\right)= \V{v}_e^{\sf{T}}\Vg{\rho}-V^{*} \leq 0\\
      & \,\,\Vg{\rho} \in  \mathcal{D}_{\text{ad}}   
\end{align}
where $g_0\left(\Vg{\rho}, \V{u}\right)$ is the objective function equal to the compliance of the mechanical system, $r\left(\V{\rho},\V{u}\right)=0$ represents the associated physical problem or the state equations written in residual form, and $g_0\left(\Vg{\rho},\V{u}\right)\leq 0$ is an additional constraint on the volume of the solid material distributed in the design domain. The vectors $\V{u}$ and $\Vg{\rho}$ represents the state solution and the material distribution respectively. As the focus here is on linear elastic problems discretized using the finite element method, the state equations $r\left(\V{\rho},\V{u}\right)=0$ can be written as 
\begin{equation}
\label{eq1002}
\M{K}\left(\Vg{\rho}\right)\V{u}=\V{f}
\end{equation}    
where $\V{f}$ is the vector with the external forces applied on the system and $\M{K}\left(\Vg{\rho}\right)$ is the so-called stiffness matrix obtained using standard finite element assembly. The computation domain is discretized using finite elements and a density value, bounded between zero and one, is assigned to each of them. Void elements are modeled by assigning the density to zero, and parts occupied with solid material are modeled with density values equal to one. All density variables, or also called design variables, are collected in the density vector $\Vg{\rho}$. The density values are allowed to vary continuously between zero and one in order to utilize gradient optimization techniques for finding a material distribution fulfilling the constraints and minimizing the objective. The vector $\V{v}_e$ collects all areas/volumes of the discrete finite elements. 

The actual physical material distribution $\Vg{\rho}_p$ is calculated by a set of transformations applied on the original density field $\Vg{\rho}$, e.g., \cite{Lazarov2016}. The first transformation is usually obtained by convoluting the density distribution with a filter functions\cite{Bendsoe2003} resulting in the so-called filtered density field $\Vg{\rho}_f$ and providing a mesh independent solution of the optimization problem. The filtered density can be utilized directly for modeling the stiffness by using the SIMP \cite{Bendsoe1989} (Solid Isotropic Material with Penalisation) material interpolation scheme where the modulus of elasticity for every element is computed as
\begin{equation}
\label{eq1002a1}
E_e=E_{\min}+\rho_e^p \left(E_{\max}-E_{\min}\right)
\end{equation}
where $E_{\max}$ is the modulus of elasticity of the solid material and $E_{\min}<<E_{\max}$ is a small regularization parameter ensuring the existence of a solution to the associated linear elastic problem, and $p$ is a penalization parameter often taken to be equal to $p=3$. 
Alternatively, as the filtered field consists of many elements with densities between zero and one, additional projection step \cite{Lazarov2016} can provide a sharper transition between void and solid. Here, the physical density is modeled directly by the filtered field, however, the presented preconditioning techniques can be applied to formulations with projections, penalization techniques different than SIMP \cite{Bendsoe2003}, to the so-called robust formulation in topology optimization \cite{Wang2011}, and other problems with high-contrast and multi-scale coefficients like level-set type of topology optimization formulations, simple parametric studies, and simulations \cite{Dijk2013}.

As stated earlier, the solution of the topology optimization problem is obtained iteratively. The optimization process starts with some admissible initial design $\Vg{\rho} \in \mathcal{D}_{\rm{ad}}$. The state solution $\V{u}$ is computed by solving \autoref{eq1002} and the gradients of the objective are evaluated by solving an adjoint equation \cite{Bendsoe2003} in the general case. For the minimum compliance problem, the gradients of the objective are given as
\begin{equation}
\label{eq1003}
\frac{\partial g_0}{\partial \Vg{\rho}_e}=-\V{u}^{\sf{T}}\frac{\partial \M{K}}{\partial \Vg{\rho}_e}\V{u}
\end{equation}
where $e$ refers to the element index. 
The densities are usually updated using the Method of Moving Asymptotes (MMA) \cite{Svanberg1987} or the optimality criteria method (OC) \cite{Bendsoe2003}. For a more detailed introduction to topology optimization, the interested readers are referred to \cite{Bendsoe2003} and \cite{Andreassen2011}.

\subsection{Physical problem and its discretization}
\label{ss:discr} 

The systems considered here (represented above as 
$r\left(\V{\rho},\V{u}\right)=0$ and simplified to
$\eqref{eq1002}$)
are linear elastic and their response is obtained by solving the Navier-Cauchy partial differential equation
\begin{align}
\nonumber
-\rm{div}\Vg{\sigma}\left(\V{u}\left(\V{x}\right)\right)& =\V{f}\left(\V{x}\right) \quad\, \V{x}\, \in \, \Omega \\
\label{eq001}
\Vg{\sigma}\left(\V{u}\left(\V{x}\right) \right)&=\M{C}\left(\V{x}\right):\Vg{\varepsilon}\left(\V{u}\left(\V{x}\right)\right)
\end{align}
where  $\Vg{\sigma}\left(\V{x}\right) $ is the stress tensor, $\Vg{\varepsilon}\left(\V{x}\right)$ is the strain tensor given by
\begin{equation}
\label{eq002}
\Vg{\varepsilon}=\frac{1}{2}\left(\nabla \V{u} +\nabla^{\sf{T}} \V{u}\right)
\end{equation}
and $\V{C}\left(\V{x}\right)$ is an elastic material properties tensor, $\V{u}\in\mathbb{R}^d$ denotes the displacement field and $\V{f}\in\mathbb{R}^d$ is the input supplied to the system, i.e., distributed and concentrated forces. The mechanical system occupies the bounded physical domain  ${ \Omega } \subset \mathbb{R}^d$.  The boundary $\Gamma={\Gamma_{D_i} \cup \Gamma_{N_i}}, i=1,\dots,d$, is decomposed into two disjoint subsets for each component $i=1,\dots,d$, $\Gamma_{D_i}$ with prescribed displacements $u_i=0$, and $\Gamma_{N_i}$ with prescribed traction $t_i$.
 
The elastic material properties tensor $\V{C}\left(\V{x}\right)$ is isotropic, and for a point $\V{x}$ in the computational domain is computed as
\begin{equation}
\label{eq003}
\M{C}\left(\V{x}\right)=E\left(\V{x}\right) \M{C}_0.
\end{equation}
In the above equation, $\M{C}_0$ is a constant tensor obtained for predefined Poisson ratio $\nu$ and modulus of elasticity one. 
For a point located in an element $K_e$ the elastic modulus is obtained using \autoref{eq1002a1}. We refer 
to \cite{Buck2013,buck2014multiscale} for design for 
domain decomposition preconditioners for the elasticity equation in the general case. 

The variational formulation \cite{Braess2007}  of \autoref{eq001}  reads 
\begin{equation}
\label{eq004}
{\rm{Find}}\,\, \V{u}\in V_0\, {\rm{s.t.}}\,\,  a\left(\V{u},\V{v}\right)=l\left(\V{v}\right)\,{\text{for all}}\, \V{v}\in V_0
\end{equation}
with bilinear form $a$ and linear form $l$
\begin{align}
\nonumber
a\left(\V{u},\V{v}\right)=\int_\Omega \left(C:\Vg{\varepsilon}\left(\V{u}\right)\right):\Vg{\varepsilon}\left(\V{v}\right) {\rm{d}} \V{x} \\
\label{eq005}
l\left(\V{v}\right)=\int_\omega\left(\V{f}\cdot\V{v}\right) {\rm{d}} \V{x} + \int_{\Gamma_N} \left(\V{t}\cdot\V{v} \right) {\rm{d}} s 
\end{align}
where $V_0=\left\{\V{v}\in\left[H^1\left(\Omega\right)\right]^d:\V{v}_i=0\,\rm{on}\,\Gamma_{D_i},i=1,\dots,d\right\}\subset V=\left[H^1\left(\Omega\right)\right]^d$. The weak formulation is discretized using finite element space $V_h \subset V_0$ with vector valued shape function defined on a uniform mesh $\mathcal{T}^h$. Each basis function is a scalar bilinear Lagrange function in one of the components and zero in the other. Substituting the shape functions in the integrals given by \autoref{eq005} for all finite elements in the mesh $\mathcal{T}^h$  leads to the linear system of equations previously introduced as \autoref{eq1002}. 

\subsection{Iterative solvers in topology optimization}
\label{ss:isolv}
Topology optimization is a computationally heavy process. The resolution of the obtained designs depends both on the transformation of the design field and the discretization of the design domain. Fine discretization is capable of representing small design features which leads to better utilization of the design freedom and at the same time to a larger system on linear equations. The total computational time is usually proportional to the number of design iterations. Every update of the design requires a negligible amount of time compared to the time necessary for solving the state problem \cite{Aage2013}. The solution of the linear system \autoref{eq1002} dominates the computational cost and requires careful selection of a solution algorithm and scalable implementation for large 3D problems \cite{Aage2017}. 

Factorization techniques are serial by nature and hard to parallelize. On the other hand, iterative linear solvers \cite{saad2003iterative} are relatively easy to parallelize and provide a scalable alternative to direct factorization techniques. The convergence rate depends on the condition number of the stiffness matrix and the clustering of the eigenvalues. The introduction of weak background material, with a modulus of elasticity orders of magnitude softer than the distributed solid, leads to an ill-conditioned system of linear equations. More precisely the condition number is of order $\eta h^{-2}$ where $h$ is a characteristic measure of the mesh size, and $\eta$ measures the contrast in the coefficients $\eta=E_{\max}/E_{\min}>>1$. Furthermore, the optimal design often consists of multiscale segments, see \autoref{figdes}, which together with the bad condition number makes the solution of the linear system extremely challenging \cite{Galvis2010, Galvis2010a, Efendiev2012a}. 

Preconditioning techniques alleviate the slow convergence speed. Therefore, the art of solving large sparse linear systems in parallel lies in the construction of computationally cheap, parallelizable and effective preconditioners. A preconditioned system $\M{M}^{-1}\M{K}\V{u}=\M{M}^{-1}\V{f}$ is obtained from \autoref{eq1002} by premultiplication with a preconditioner $\M{M}$. The most effective preconditioner $\M{M}^{-1}$ is the exact inverse of the stiffness matrix $\M{K}$. However, direct construction requires matrix factorization which as already discussed is not scalable and is extremely expensive for large problems. Thus, a preconditioner in the form of a multigrid procedure \cite{Vassilevski2008} or a domain decomposition \cite{Toselli2004} provides the most efficient solution procedure. 

\subsection{High-contrast coefficients in topology optimization}
\label{ss:hctop}
The focus here is on domain decomposition preconditioners, in particular, on two-levels domain decomposition techniques.  The classical variants of these preconditioners \cite{Toselli2004} do not perform well for high-contrast problems \cite{Efendiev2012a}. The condition number estimates for the traditional domain decomposition case depends on the contrast $\eta$ if the high-conductivity regions are not aligned with the coarse mesh of the decomposition. We refer to \cite{Efendiev2012a,Efendiev2013a,Buck2013,buck2014multiscale} where it is demonstrated  that if the material properties with local variations are enclosed in a coarse block the performance of the standard preconditioner is not affected by the contrast. However, for cases with several extended high-conductivity areas (high stiffness regions for linear elasticity) crossing the coarse block boundaries, the performance deteriorates significantly. That is precisely the case for topology optimized linear elastic structures. These design features can be observed in \autoref{figdes} as well as in previous articles on the topic \cite{Lazarov2014, Alexandersen2015a, Alexandersen2015}.

We apply the Generalized Multiscale Finite Element methods (GMsFEM) framework introduced in~\cite{Efendiev2011a, Galvis2010a, Efendiev2013a} to construct robust and fast solution algorithms adapted to linear elasticity.  Similar to the constructions of diffusion type multilevel domain decomposition preconditioners, the proposed design depends on the behavior of the coefficient inside local coarse node neighborhoods. We demonstrate that the preconditioned solvers perform, in terms of iterations and condition number estimates, independently of the contrast in the media properties.

The most important part of the construction of a multilevel domain decomposition preconditioner  is well know to be the coarse level. The coarse level of the preconditioner has to provide a good local approximation of the kernel of the elasticity operator \cite{Efendiev2010}. Also, it should contain all eigenmodes with corresponding eigenvalues dependant on the contrast of the coefficient.  Thus, an adequately chosen eigenvalue problem is constructed and solved locally to ensure the desired behavior. 

Two major approaches can be identified in the current literature. For coarse spaces with standard dimension (for scalar problems that is, one basis function per coarse node and diffusion coefficient $\kappa$) it is imperative to have coarse basis functions $\phi$ such that $\kappa |\epsilon (\phi)|^2$ is bounded independently of the contrast. For high-conductivity regions (high stiffness regions for linear elasticity)  restricted within the coarse element, i.e.,  there are not any long channels crossing the edges of the coarse element, the above requirement is fulfilled by classical multiscale finite element basis functions obtained by energy minimization \cite{Buck2013}.  For high-contrast coefficients with extended channels, even for isotropic problems, the above condition cannot be fulfilled. Therefore, to achieve robust behavior, an enrichment procedure is implemented by adding basis functions that approximate the contrast dependent eigenmodes of the operator locally.

Apart from the fact that the material coefficients in topology optimization show multiscale variations combined with high-contrast, an additional complication comes from the fact that throughout the optimization iteration the density field and correspondingly $E(x)$ are changing as the optimization converges to the optimized design. A topology represented by a particular density distribution evolves with the iterations leading to iteration dependent multiscale structure. High-contrast channels may break apart or joint together during the optimization iterations within a globally connected subdomain $\Omega^{mat}$ restricted to a coarse neighborhood. Dealing with such an additional complication requires re-computation of the preconditioner as the optimization iteration advances towards the final solution. Similar to \cite{Lazarov2014, Alexandersen2015}, we pay extra attention to the building cost of the preconditioner, especially to the construction of the coarse level. Therefore, two new alternatives are proposed to reduce the cost of recomputing the basis functions for the coarse space:

\begin{itemize}

\item{Computation of eigenvalue problems using a randomized algorithm. The randomized approximation of the local eigenvalue problems for GMsFEM is proposed in~\cite{Calo2016b, Efendiev2013a} and realizes a significant reduction of the cost of computing the coarse basis functions that generate the coarse space. A detailed introduction can be found in the overview article \cite{Halko2011}.}

\item{ The utilization of a preconditioner constructed for the heat equation in a similar high-contrast multiscale media to precondition the elasticity equation. Here, we follow the ideas presented in \cite{Gustafsson1998, Gustafsson2002} for homogeneous media. More precisely, we utilize the coarse spaces generated for the heat equation with a diffusion coefficient defined using the solid material region $\Omega^{mat}$. The local eigenvalue problem is related to the diffusion operator and not to the elasticity equation operator. Therefore, the local eigenvalue problem, for the same resolution, is twice smaller in two-dimensions and three times smaller in three dimensions reducing the computational cost significantly. The combination with randomized algorithm results in a solver for the topology optimization an order of magnitude faster than the one presented earlier in \cite{Alexandersen2015}.}
\end{itemize}

\begin{remark}
An unexpected advantage of using the heat equation to precondition the elasticity equation in the context of GMsFEM framework is that in the case of the heat operator it is easier to define threshold strategies for the decay of the local eigenvalues, and therefore adaptivity to coefficients is easier to obtain.
\end{remark}



\subsection{Topology optimization example}

\begin{figure}
\center{
\includegraphics[height=0.43\textwidth]{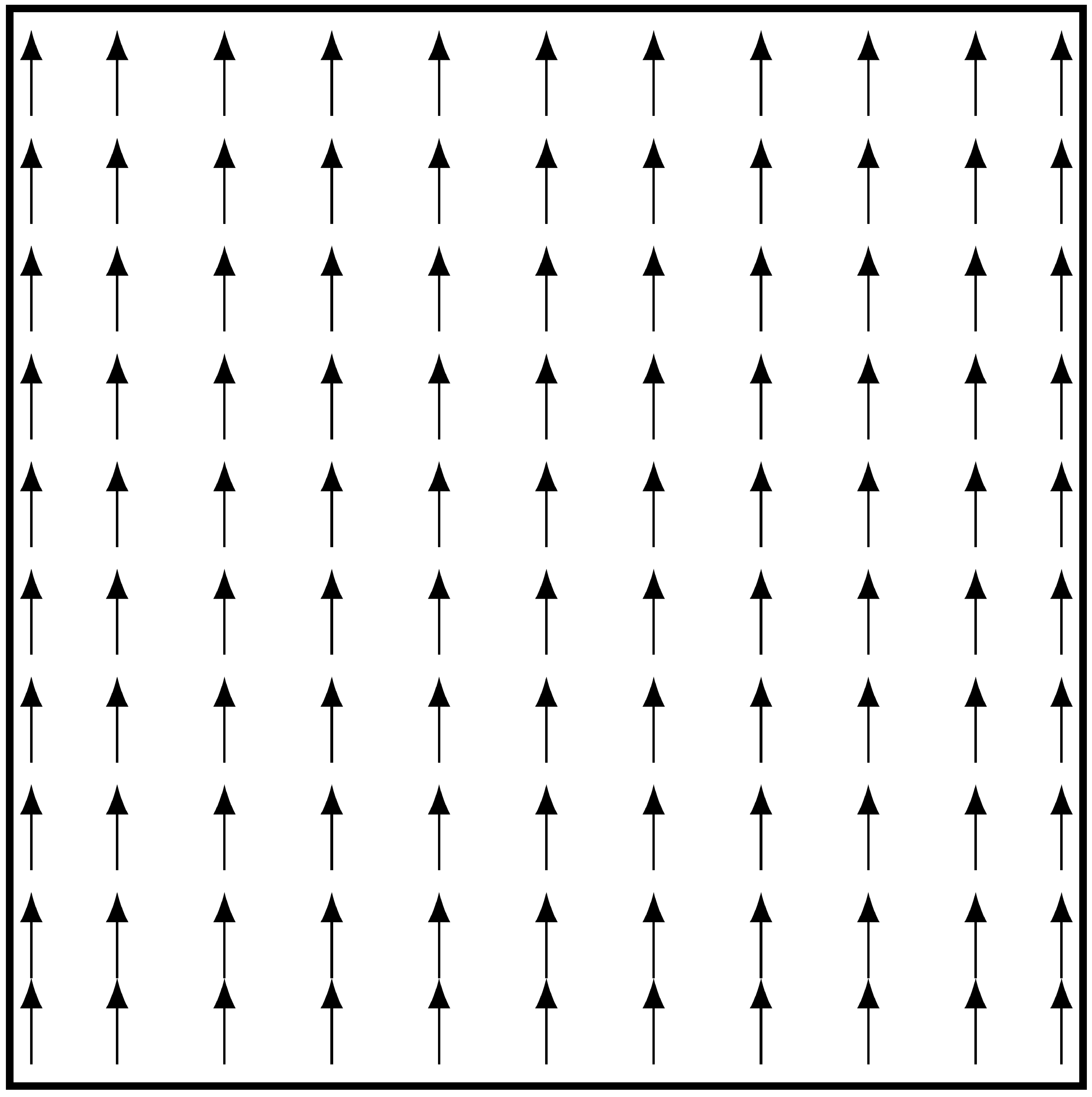}\quad
\includegraphics[height=0.43\textwidth]{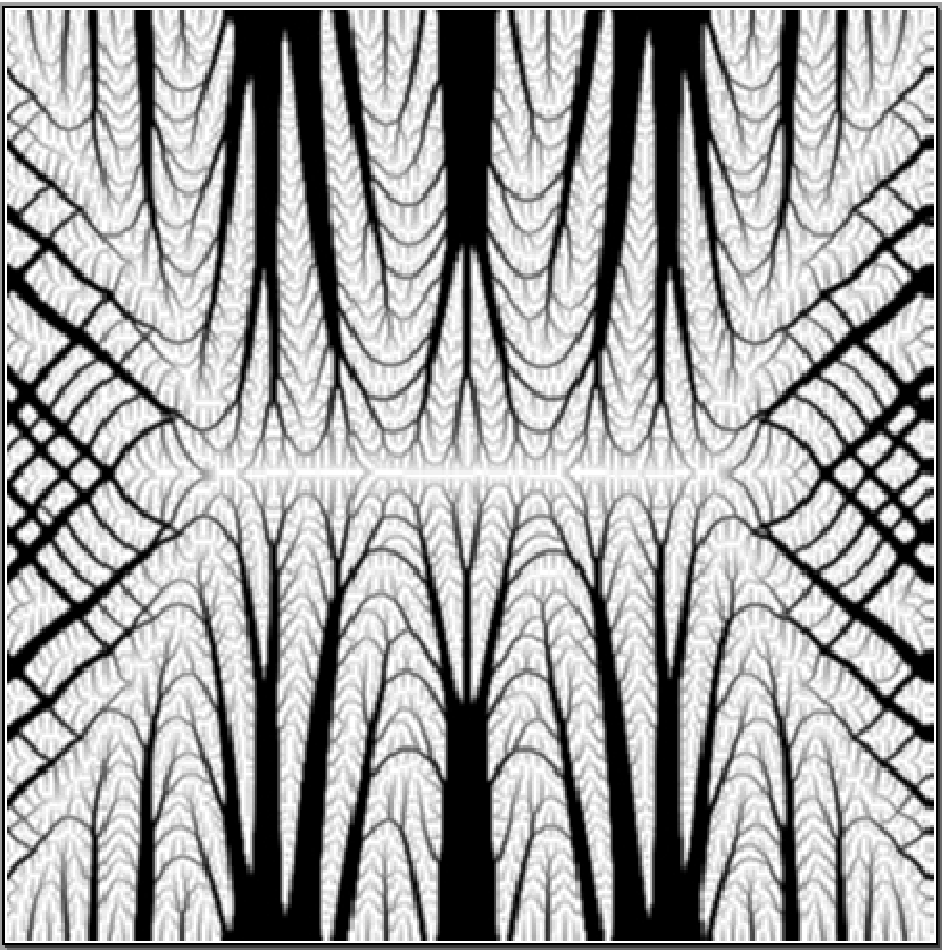}}
\caption{Design domain and topology optimization of a 2D plane stress problem with homogeneous Dirichlet boundary conditions and distributed force over the domain.}
\label{figdes}
\end{figure}

To demonstrate the performance of the preconditioner we consider a square 2D plane stress example with homogeneous Dirichlet boundary condition and distributed force over the domain, as shown in \autoref{figdes}. The setup even though difficult to find in the engineering practice is an excellent example to test preconditioners as the optimization results in a highly sophisticated multiscale design. The design domain is partitioned on $400\times 400$ square mesh elements with $20\times 20$ coarse elements. The maximum number of eigenvectors per coarse neighborhood is taken to be $6$. The topology optimization algorithm was set to stop at $300$ iterations. 

As the focus here is on the effective preconditioning techniques, we limit the topology optimization problems to the one presented in \autoref{figdes}. The presented topology is obtained using the standard and the preconditioned solvers presented here. Depending on the preconditioner update strategy, the number of the iterations reduced an order of magnitude. The fastest total solution time is obtained by employing the strategy presented in \cite{Alexandersen2015}, i.e., either updating the preconditioner after a fixed number of optimization iterations or updating it after the number of the iterations of the local iterative solver exceeds a predefined threshold. More extensive 3D topology optimization studies will be presented in the following articles, and here we will focus on demonstrating the numerical properties of the developed preconditioners and the associated computational complexity.

The proposed replacement of a standard elasticity with a diffusion solver for the local eigenvalue problem results in a time reduction factor of $3^3$. The above result emphasizes the contribution of the local solves to the total solution time and is in line with the theoretical predictions. The diffusion problem is three times smaller compared to the elasticity, and the cost of the local eigenvalue solves is proportional to $n^3$, where $n$ is the number of local degrees of freedom. The randomized solver reduces further the computational cost resulting in an additional factor of two to three. It should be pointed out that the problems considered here are 2D, and the randomized modification is expected to have a more significant effect on the preconditioner time in larger 3D problems. 

\section{GMsFEM two-levels domain decomposition for the elasticity equation}

The focus In this section is on the utilization of GMsFEM coarse spaces in the construction of two-levels domain decomposition preconditioners. In particular, we show that the proposed preconditioners yield a contrast-independent condition number, and thus they are optimal in terms of physical parameters. Additional theoretical details can be found in \cite{Galvis2010, Efendiev2011} and further extensions of the results to multilevel methods are presented in \cite{Efendiev2011c}. 

\begin{figure}
\centering
\begin{tikzpicture}[scale=1.2]
    \filldraw [draw=black, fill=gray!25!white, very thick] (0,0) rectangle (5,5);
    
    \draw [step=0.5, black!70!white, dashed] (0,0) grid (5,5);
    \draw [black, very thick] (0,0) rectangle (5,5);
    
    \filldraw [draw=black, fill=gray!25!white, very thick] (7,1) rectangle (10,4);
    \draw [step=0.5, black!70!white, dashed] (7,1) grid (10,4);
    \draw [black, very thick] (7,1) rectangle (10,4);
    \fill [gray] (3.5,1) rectangle (4,1.5);
    \draw[-Latex, thick] (4,1.5) to [bend left] (7,2.5);
    \draw[|-|](3.5,0.9)--(4,0.9);
    
    \draw(3.75,1.25) node {$D_i$};
    \draw(3.75,0.9) node [anchor=north]{$H$};
    
    \draw[|-|](7,0.85)--(10,0.85);
    \draw(8.5,0.85) node [anchor=north]{$H$};
    
    \draw[|-|](6.85,1.5)--(6.85,2);
    \draw(6.85,1.75) node [anchor=east]{$h$};
\end{tikzpicture}
\caption{Decomposition of the computational domain into non-overlapping coarse blocks (subdomains) $D_i$.}
\label{fig:dom01}
\end{figure}
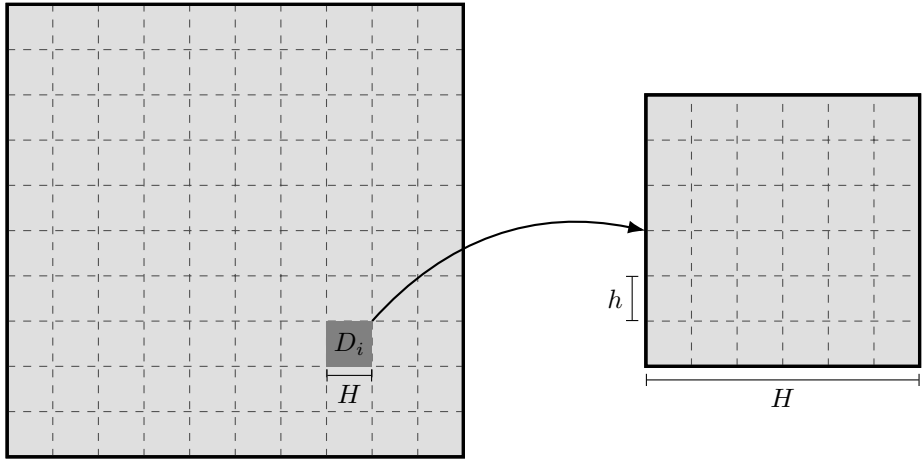

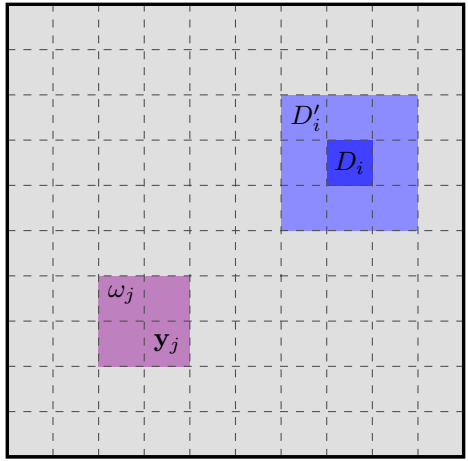
\begin{figure}
\centering
\begin{tikzpicture}[scale=1.2]
    \filldraw [draw=black, fill=gray!25!white, very thick] (0,0) rectangle (5,5);
	\fill [violet!50!white] (1,1) rectangle (2,2);

    \fill [blue!45!white] (3,2.5) rectangle (4.5,4);
    \fill [blue!75!white] (3.5,3) rectangle (4,3.5);

    \draw(3,4) node [anchor=north west]{$D^\prime_i$};
    \draw(3.75,3.25) node {$D_i$};
    \draw [step=0.5, black!70!white, dashed] (0,0) grid (5,5);
    \draw [black, very thick] (0,0) rectangle (5,5);
    \node[mark size=1.5pt, color=black] at (1.5,1.5) {\pgfuseplotmark{*}};
    \draw(1.5,1.25) node [anchor=west]{$\mathbf{y}_j$};
    \draw(1,2) node [anchor=north west]{$\omega_j$};
\end{tikzpicture}
\caption{Definition of overlapping domain decomposition and a neighborhood $\omega_j$ of a coarse node $\V{y}_j$.}
\label{fig:dom02}
\end{figure}

The computational domain, shown in \autoref{fig:dom01}, is partitioned into finite elements which are agglomerated into larger non-overlapping blocks $D_i$. Based on the above decompotions an overlapping decomposition, shown in \autoref{fig:dom02}, denoted as $\{D_i'\}_{i=1}^N$ is obtained from an original non-overlapping decomposition $\{ D_i\}_{i=1}^N$ by enlarging each subdomain $D_i$ to 
\begin{equation}
\label{eq:def:overlapping-subdomains}
D'_i=D_i\cup \{x\in D, \mbox{dist}(x,D_i)<\delta_i\}, \quad i=1,\dots,N,
\end{equation}
where $\mbox{dist}\left(\cdot\right)$ is some distance function. The overlapping subdomains $\{D_i'\}$ and the coarse triangulation $\mathcal{T}^H$ are not related in general settings. Two partitions of unity covering the whole computational domain, one for $\{D_i'\}$ and one for $\mathcal{T}^H$, can be chosen independently of each other. However, in the numerical experiments presented later in the paper, we assume that the overlapping subdomains $\{D_i'\}$ coincide with the coarse vertex neighborhoods $\{\omega_i\}$ of $\mathcal{T}^H$, and in this case $\delta \asymp H$, where $\delta=\max_{1\leq i\leq N} \delta_i$ is the overlapping parameter.

Based on the above decomposition, the preconditioned operator is $\M{M}_{EE}^{-1}\M{K}$, and the preconditioner matrix is defined as 
\begin{equation}\label{eq:2lEp}
 \M{M}_{EE}^{-1}= \M{M}^{-1}_{E,1}+ \M{M}^{-1}_{E,2},
\end{equation}
The part corresponding to the first level is
\begin{equation}\label{eqjuan:localproblemsE}
\M{M}^{-1}_{E,1}r= \sum_{i=1}^N \M{R}_{i}^{\top} \M{K}_{i}^{-1}\M{R}_{i} \V{r},
\end{equation}
where $\M{K}_{i}= \M{R}_{i}\M{K}\M{R}_i^T$, 
$1\leq i\leq N$, and the part corresponding to the second (or coarse) level is
\begin{equation*}
\M{M}^{-1}_{E,2}r= \M{R}_{0}^{\top} \M{K}_{0}^{-1}\M{R}_{0}\V{r}
\end{equation*}
where $\M{K}_{0}= \M{R}_{0}\M{K}\M{R}_0^T$. The matrix $\M{R}_0^T: V_0\rightarrow V$ consists of vectors defining the so-called coarse space $V_0$ which is obtained by interpolating the coarse functions onto the fine mesh. The matrices $\M{R}_i^T: V_i\rightarrow V$ are rectangular and consist of zeros and ones and utilized to extract the degrees of freedom that lie inside the subdomains $D'_i$. 


The fine-scale linear system \autoref{eq1002} is solved iteratively with the preconditioned conjugate gradient (PCG) method. The application of the preconditioner  involves solving a coarse-scale system $\M{M}^{-1}_{E,2}\V{r}$ and solving a set of local problems $\M{M}^{-1}_{E,1}\V{r}$  in each iteration. The main goal is to reduce the number of iterations in the solution process. Without the coarse space $\M{R}_0$, the preconditioner usually acts as a smoother and the convergence depends on the number of the coarse blocks. Thus, the appropriate construction of the coarse space $V_0$ plays a crucial role in obtaining robust iterative domain decomposition methods and ensures iteration number independent on the contrast and the characteristic mesh size. See \cite{Toselli2004, Efendiev2012a,Efendiev2013a}


\subsection{Generalized multiscale coarse finite element spaces for the elasticity equation}\label{juansec:gmsfemE}

The construction of the coarse space $V_0\subset V$ for linear elasticity follows the strategy outlined in \cite{Efendiev2011, Efendiev2013c} for the diffusion equation. The process starts with the selection of an initial set of basis functions ${\chi}_{i}$ that form a partition of unity and are associated with the coarse domains $\omega_i$. Additional sets of coarse basis functions $\psi_{i,l}$ are defined with respect to the fine mesh $\mathcal{T}_h$. These are computed by solving a local eigenvalue problem
\begin{equation}\label{eq:juanelasticityeigenproblem}
-\rm{div}\Vg{\sigma}\left( \Vg{\psi}_{i,l} \left(\V{x}\right) \right) =  \lambda_l E\left(\V{x}\right) \Vg{\psi}_{i,l}\left(\V{x}\right), \quad \,\,\, \V{x}\, \in  \omega_l 
\end{equation}
with homogeneous Neumann boundary condition on $\partial \omega_l$ and Dirichlet boundary condition on $\partial \omega_l\cap \partial D$ if $\partial \omega_l\cap \partial D$ is not empty. In matrix form we have
\begin{equation}
\label{eq3001}
\M{K}_E^{l}\Vg{\phi}=\lambda \M{M}^{l} \Vg{\phi}
\end{equation}
with eigenvalues and the eigenvectors denoted as $\left\{ \lambda_i^{\omega_l}\right\}$ and $\left\{\Vg{\phi}_i^{\omega_l}\right\}$ respectively. The eigenvalues are ordered as
\begin{equation}
\lambda_1^{\omega_l}\leq\lambda_1^{\omega_l}\leq\lambda_1^{\omega_l}\leq\dots\leq\lambda_j^{\omega_l}\leq\dots
\end{equation}
The multiscale coarse basis functions are constructed by multiplying selected eigenfunctions $\left\{\psi_j^{\omega_l}\right\}$ with the partition of unity function $\chi_l$ associated with subdomain $\omega_l$. The coarse GMsFEM space for the entire problem is defined by 
\begin{equation}\label{juaneq:VOE}
V_{0,E}=\rm{span}\left\{\Psi_{i,l}=\chi_l\Vg{\psi}_{i,l}^{\omega_l},\,\, i=1,\dots,N,\,\, l=1,\dots,L \right\}
\end{equation}
where $N$ and $L$ denote the number of eigenvectors and the number of coarse blocks $\omega_l$ respectively. 

\begin{figure}
\centering
\includegraphics[width=1\textwidth]{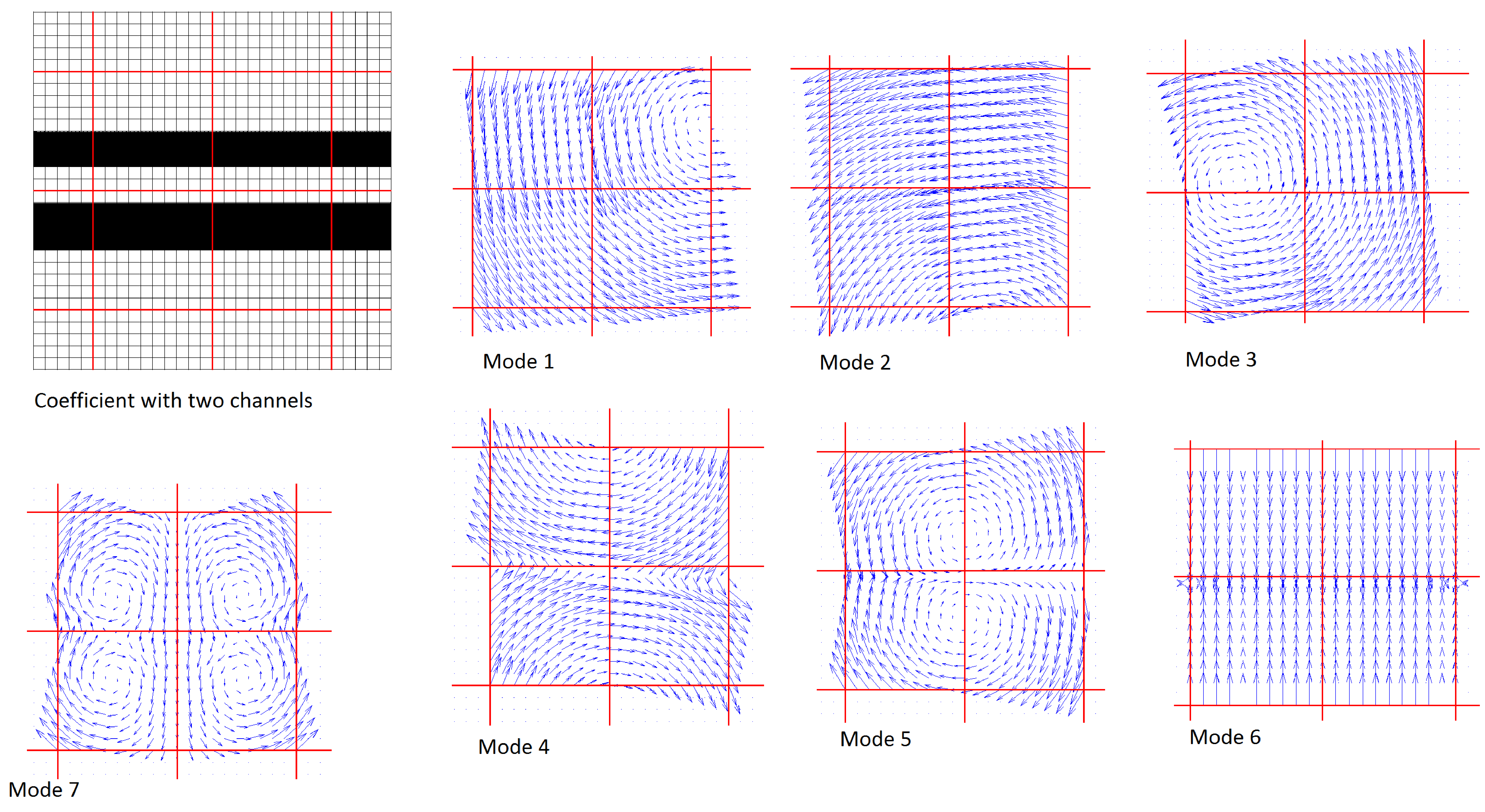}
\caption{Two disconnected high-stiffness regions and their contrast dependent modes. Three for each region corresponding to the RBMs. The next mode is contrast independent as observed in our numerical tests. Note that, restricted to high-contrast regions, these contrast dependent modes correspond to linearly independent rigid body motions.}
\label{fig:elasticitymodes}
\end{figure}

For the diffusion case, the contrast dependent eigenvalues are well separated, i.e., a clear jump can be observed between the contrast dependent and the rest of the eigenmodes computed on a given patch $\omega_l$. Thus, selecting the low order contrast dependent modes results in a preconditioner which provides a condition number for the preconditioned operator independent of the contrast. However, for the linear elastic case, such behavior for the low order modes is difficult to be observed. The optimal number of low order modes is related to the disconnected high-stiffness regions and the RBM (rigid body motion) of the region \cite{Efendiev2012a}. We illustrate this fact in \autoref{fig:elasticitymodes} where we picture a coarse node neighborhood with two disconnected high-stiffness regions and its contrast dependent modes. There are three modes for each high-stiffness region, and these correspond to the rigid body modes. The next mode, the one corresponding to the eigenvalue number seven in increasing order, is contrast independent as observed in our numerical tests. Note that, restricted to high-contrast regions, the first 6 modes correspond to three linearly independent rigid body modes per inclusion. The space RBM of rigid body modes on a set $\Omega \subset \mathbb{R}^d$ is defined for $d=2$, by
\begin{equation*}
    \rm{RBM}(\Omega)=\left\lbrace v\in [L^2(\Omega)]^2: v=a+b\left( \begin{array}{c}
         -x_2  \\
         x_1
    \end{array}\right), a\in \mathbb{R}^2, b\in \mathbb{R}, x\in\Omega \right\rbrace.
\end{equation*}
\begin{remark}
In most of our numerical experiments using the eigenvalue problems in \autoref{eq3001} we did not find any automatic way to implement a threshold to select the adequate number of modes in each coarse neighborhood. The numerical experiments were performed by specifying the number of basis functions based on the number of disconnected high-stiffness regions present in the specific coarse node neighborhood. Later we will introduce the use of GMsFEM basis functions constructed for the heat equations where an automatic threshold can be implemented to select the adequate number of basis functions in each coarse neighborhood. 
\end{remark}

\begin{figure}[H]
\begin{minipage}{0.49\textwidth}
\centering
\includegraphics[width=0.6\textwidth]{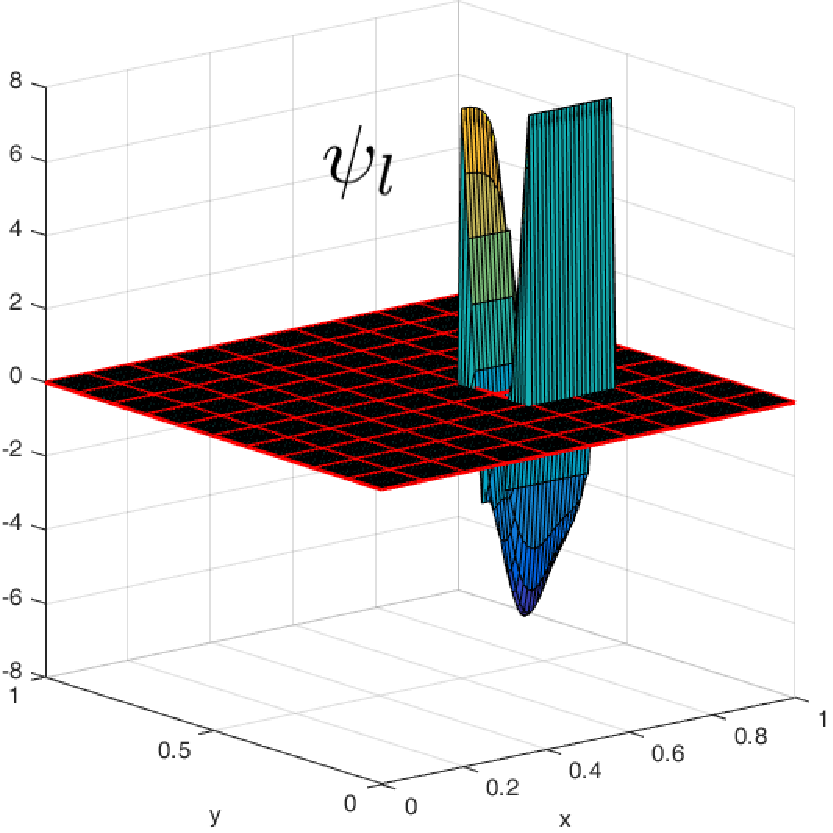}
\end{minipage}
\begin{minipage}{0.49\textwidth}
\centering
\includegraphics[width=0.6\textwidth]{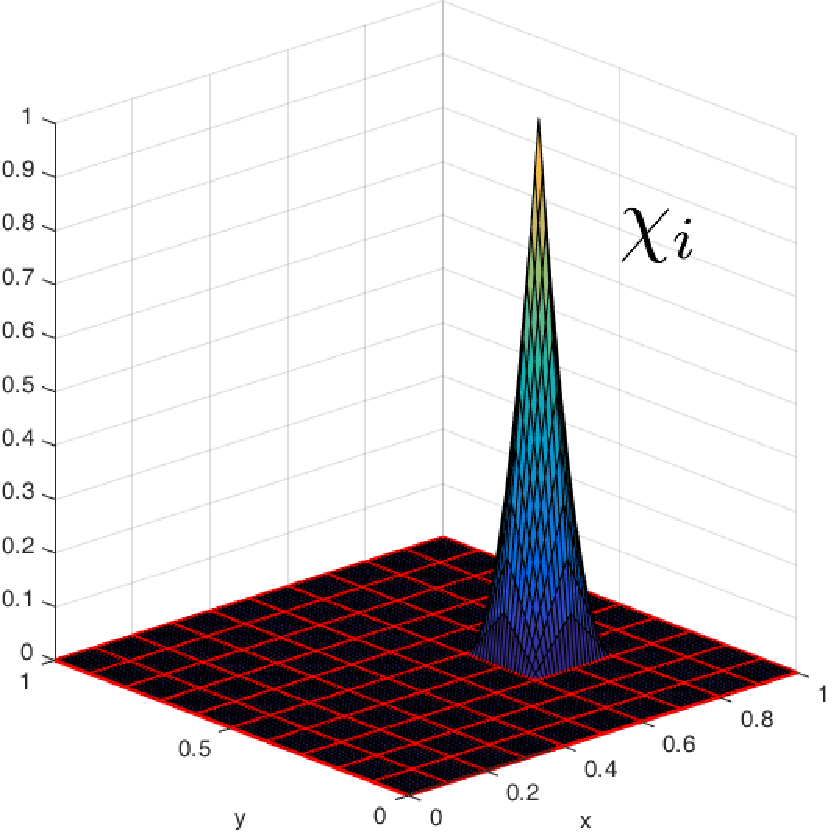}
\end{minipage}
\centering
\includegraphics[width=0.3\textwidth]{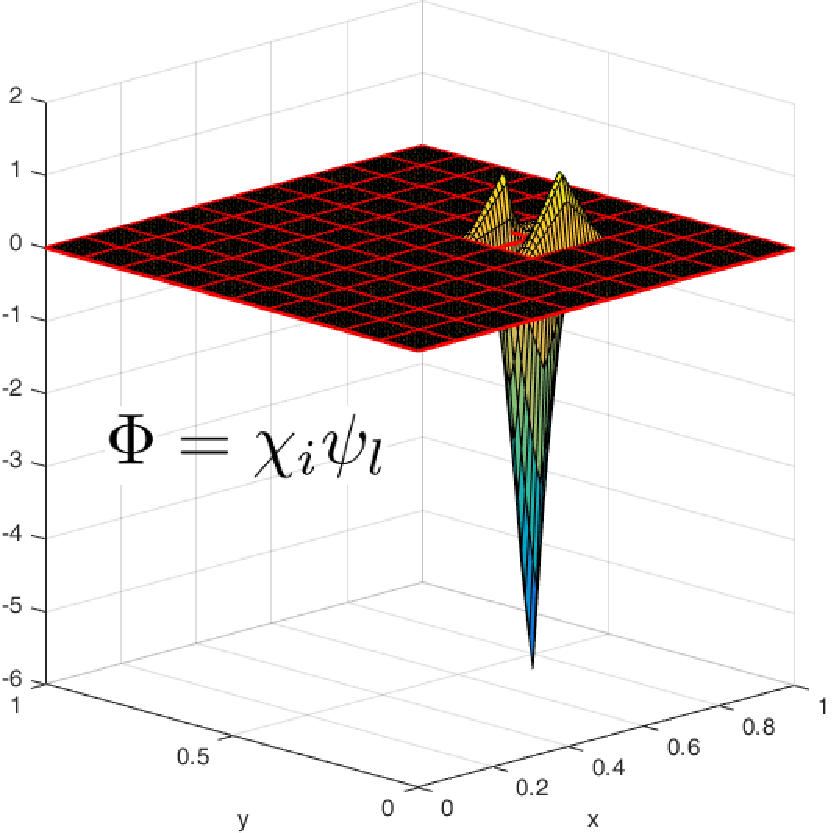}
\caption{Basis function construction for an element of the domain. To illustrate the process only one displacement component, along the x-direction, of the eigenvector is shown in the coarse region. The values for the other components are obtained following the same procedure.}\label{fig:basisillustration}
\end{figure}

For the elasticity problem the GMsFEM approximate the solution on the coarse grid as,
\begin{equation*}
\V{u}_0=\sum_{i,\ell }^{N} c_{i,\ell} \Vg{\Phi}_{i,\ell}
\end{equation*}
where $c_i$ are the unknown constants. The coarse matrix is constructed, 
\begin{equation*}
\M{K}_{0}= \M{R}_{0}\M{K} \M{R}_{0}^{\top},
\end{equation*}
where 
\begin{equation*}
\V{R}_{0}^{\top}=[\Vg{\Phi}_1 \quad \Vg{\Phi}_2 \quad \cdots \quad \Vg{\Phi}_{N_c}].
\end{equation*}
and the multiscale finite element solution is the finite element projection of the fine-scale solution on $V_{0,E}$, that is, $\M{R}_{0}\V{u}_0$ where
\begin{equation*}
\M{K}_{0}\V{u}_0=\V{f}_0,
\end{equation*}
and $\V{f}_0=\M{R}_{0}^{\top}\V{f}$.
\section{A randomized algorithm for eigenvalues computation}
\label{rndalg}
Inspired by \cite{Halko2011}, randomized algorithms have been introduced in GMsFEM in \cite{Efendiev2013c,Calo2016b}. See also 
\cite{galvis2017overlapping}. The basic idea is to utilize random sampling to generate low-rank approximations to the set of matrices utilized for finding the coarse shape functions \autoref{eq3001}. The main difference between deterministic factorization techniques is the convergence criteria. For randomized algorithms, the converge is probabilistic. However, as stated in \cite{Halko2011}, the probability of failure is often negligible of order $10^{-15}$. Therefore, for practical applications involving singular values or eigenvalues decompositions, a randomized algorithm offer a computationally cheap alternative to direct factorization methods.

The algorithm for finding approximations of the eigenvalues and the eigenvectors starts with dimension reduction. For each subdomain $\omega_i, i=1,\dots,N_s$ the following sequence of steps is executed in parallel

\begin{itemize}

\item{Generate forcing terms $\V{f}_1,\V{f}_2,\dots,\V{f}_M$ using (for instance) an uniform random distribution and $\int_{\omega_i} \V{f}\left(\V{x}\right){\text d} \V{x} =\Vg{0}$.}

\item{Compute local solutions $\M{K}_E^i \V{u}_l=\V{f}_l, l=1,\dots,M$}

\item{Generate $\mathcal{W}_i=\left\{\text{span}\left\{{\V{u}_l}\right\}\cup \text{RBM}\left(\omega_i\right)\right\}$}

\item{Solve a reduced eigenvalue problem $\widetilde{\M{K}}_E^i \widetilde{\Vg{\phi}}= \widetilde{\lambda}\widetilde{\M{M}}_E^i \widetilde{\Vg{\phi}}$}

\end{itemize}

The reduced size matrices are generated as:
\begin{align}
\widetilde{\M{K}}_E^i &=\M{U}_i^{\sf{T}} \M{K}_E^i \M{U}_i\\
\widetilde{\M{M}}_E^i &=\M{U}_i^{\sf{T}} \M{M}_E^i \M{U}_i
\end{align}
where every column of $\M{U}_i$ is a vector from $\mathcal{W}_i$. The approximations of the eigenvalues and the eigenvectors are computed as 
\begin{align}
\lambda^{\omega_i} &=\widetilde{\lambda^{\omega_i}}\\
\Vg{\phi}^{\omega_i} &=\M{U}_i \widetilde{\Vg{\phi}^{\omega_i}}.
\end{align} 
Usually, the matrix $\M{U}_i$ holds several vectors obtained by singular value decomposition (SVD) of the snapshot matrix $\left[\V{u}_1,\V{u}_2,\dots,\V{u}_M\right]$ enriched with the three rigid body modes in 2D and six in 3D. The actual number of vectors depends on the desired number of shape functions. As discussed in \cite{Halko2011}, see also \cite{Calo2016b}, for a target number of shape functions $k$, the rank of $\M{U}_i$ can be selected to be as low as $k+5$. The computational time of a standard generalized eigenvalue algorithm is proportional to $n^{3}$, where $n$ is the size of the problem. Thus, the reduced problem offers significant speed up. However, the cost associated with the solution of the linear system cannot be reduced further. A possible speed up based on splitting of the displacement field is discussed in the following section. 

\section{Displacement field splitting preconditioner}
\label{dispsp}
The basic idea behind the proposed preconditioner for linear elastic isotropic problems is presented in \cite{Axelsson1978, Blaheta1994, Gustafsson1998, Gustafsson2002}. The displacement vector is split in two blocks (three blocks in 3D) using the displacement components aligned with x, and y-directions, i.e., $\V{u}=\left[\V{u}_x, \V{u}_y \right]^{\sf{T}}$. Following the above decomposition, the stiffness matrix and the load vector $\V{f}$, from \autoref{eq1002},  can be written in a block form as 
\begin{equation}
\label{eq4001}
\left[\begin{array}{cc} \M{K}_{xx}& \M{K}_{xy}\\
\M{K}_{xy}^{\sf{T}}& \M{K}_{yy}  \end{array} \right]\left[ \begin{array}{c} \V{u}_x\\ \V{u}_y \end{array}\right]=\left[\begin{array}{c} \V{f}_x\\ \V{f}_y \end{array}  \right].
\end{equation}
The displacement split preconditioner matrix is constructed by keeping only the diagonal block matrices.
\begin{equation}
\label{eq4002}
\M{C}_{\text{EL}}=\left[\begin{array}{cc} \M{K}_{xx}& \Vg{0}\\
\Vg{0}& \M{K}_{yy} \end{array} \right].
\end{equation}
As demonstrated in \cite{Gustafsson1998, Gustafsson2002}, the condition number of the preconditioned operator depends only on the Poisson's ratio and is given as
\begin{equation}
\label{eq4003}
\kappa\left(\M{C}_{\text{EL}}^{-1}\M{K} \right)=\frac{2}{1-\widetilde{\nu}}
\end{equation}
where $\widetilde{\nu}=\frac{\nu}{1-\nu}$. Furthermore, the condition number does not depend on the contrast of the materials distributed in the computational domain. The latter makes the block diagonal preconditioner perfect for topology optimization problems. Close inspection of the diagonal block matrices reveals that they are equal for isotropic elastic problems $\M{K}_{xx}=\M{K}_{yy}$, and are equivalent to a stiffness matrix obtained by the discretization of a scalar Laplace problem. The last property is utilized in the proposed preconditioners to reduce further the computational cost. Instead of constructing the coarse bases by solving the linear elastic eigenvalue problem given by \autoref{eq3001}, the idea here is to solve a reduced eigenvalue problem associated with the Laplace problem.

The MsFEM preconditioner and the procedure for finding the coarse bases can be applied directly to the diagonal blocks, can be combined with randomized algorithms, or can be utilized to form a coarse space and subsequently project the full elastic matrix. The above options lead to different preconditioners discussed in detail in the following and numerically tested in \autoref{sec::num}.
\subsection{Local problems and coarse basis from a diffusion operator}
The first preconditioner derived from the field splitting preconditioner consists of a fine level block-diagonal preconditioner constructed for the diffusion equation and a coarse-level part obtained from a coarse set of basis functions obtained again using an eigenvalue problem for the diffusion case. The preconditioner can be written as 
\begin{equation}
\M{M}_{HH}^{-1}=\left[\begin{array}{cc} \M{M}_H^{-1} & \Vg{0} \\ \Vg{0}& \M{M}_H^{-1} \end{array}\right] + \left[\begin{array}{c} \M{R}_{H}  \\ \M{R}_H \end{array}\right] \M{K}_{E,0}^{-1}  \left[\begin{array}{c} \M{R}_{H} \\ \M{R}_H \end{array}\right]^{\sf{T}}
\label{eq2025}
\end{equation}
where
\begin{equation}\label{eq:juan001c}
\M{K}_{E,0}=\left[\begin{array}{c} \M{R}_{H} \\ \M{R}_H \end{array}\right] \M{K}_E \left[\begin{array}{c} \M{R}_{H} \\ \M{R}_H \end{array}\right]^{\sf{T}}
\end{equation}
and $\M{K}_E$ is the stiffness matrix obtained by a finite element discretization of the Navier-Stokes equations for linear elasticity. The index $H$ denotes quantities obtained from the diffusion case and index  $E$ quantities obtained from the elasticity equations. For instance, the image of the operator (or column space of the  matrix) $R_H^T$ is the coarse space $V_{0,H}$ constructed using the GMsFEM procedure but starting with the diffusion operator $-\mbox{div}(E(x)\nabla (\cdot))$. In this case the construction of the coarse basis function is done component-wise (x-and y-direction displacements) and each of those uses the local eigenvalue problem 
\begin{equation}\label{eq:juaneheateigenproblem}
-\rm{div}(\kappa(x)\left(\psi_{i,l} \left(\V{x}\right) \right) =  \lambda_l \kappa(x)\left(\V{x}\right) {\psi}_{i,l}\,\,\, \V{x}\, \in  \omega_l 
\end{equation}
where $\kappa=\tr(E)$ or some other quantity that captures the high-contrast and multiscale structure of $E$. This eigenvalue problems is posed with homogeneous Neumann boundary condition on $\partial \omega_l$ and Dirichlet boundary condition on $\partial \omega_l\cap \partial D$ if $\partial \omega_l\cap \partial D$ is not empty. In this case we have 
\begin{equation}\label{juaneq:V0H}
V_{0,H}=\rm{span}\left\{\chi_l[\psi_{i,l}^{\omega_l},0], 
\chi_l[0,\psi_{i,l}^{\omega_l}],\,\, i=1,\dots,N,\,\, l=1,\dots,L \right\}.
\end{equation}
See Section  \ref{juansec:gmsfemE} for comparison with the construction of $V_{0,E}$ that uses local elasticity eigenproblems. The cost of constructing 
$V_{0,H}$ is less than that of constructing $V_{0,E}$ in 
\eqref{juaneq:VOE} since we solve smaller local eigenvalue problems. Additionally, we also observed numerically that the selection of contrast-dependent modes can be performed automatically in the case of the local heat operator eigenvalue problem which is harder to do in the case of the local elasticity operator eigenproblem. For the numerical test, we consider the approximation of the eigenvalue 
problem \eqref{eq:juaneheateigenproblem} using a randomized method 
similar to the one described in Subsection \ref{rndalg}.

The first level of the preconditioner 
\eqref{eq2025} contains the operator 
$\M{M}_H^{-1}$ that, in a similar manner to that of the construction of the coarse basis, it denotes the additive first level of the precoditioner constructed for the diffusion operator 
$-\mbox{div}(E(x)\nabla (\cdot))$. Applying 
$\M{M}_H^{-1}$ requires the solution of Dirichlet local  diffusion equations and then add their extensions by zero outside the corresponding subdomain.

\subsection{Local problems and coarse basis from a diffusion operator with rigid body motions enrichment}

The previous preconditioner includes  the null space for the Laplace operator, which is capable of representing only the rigid body translations. The null space of the elasticity operator is larger. The coarse space must be able to represent all rigid body modes to ensure convergence independent on the problem size \cite{Griebel2003}. Thus, the coarse space is enriched with additional vector in 2D for the rigid body rotations. The construction of the vector elements is based on the vector $\V{r}=\left[\V{r}_x, \V{r}_y\right]^{\sf{T}}$ which is the interpolation on the fine-grid of the vector function representing the rotation. Thus, to the basis functions previously constructed we add the vector $\left[\chi_i \V{r}_x, \chi_i \V{r}_y \right]$ obtained for each coarse neighborhood.  For 3D problems, the coarse basis should include three additional rigid body rotational modes. The modified preconditioner, in this case, is given as 

\begin{equation}
\M{M}_{HH+Rot}^{-1}=\left[\begin{array}{cc} \M{M}_H^{-1} & \Vg{0} \\ \Vg{0}& \M{M}_H^{-1} \end{array}\right] + \left[\begin{array}{c} \M{R}_{H+Rot}  \\ \M{R}_{H+Rot} \end{array}\right] \M{K}_{E,0}^{-1}  \left[\begin{array}{c} \M{R}_{H+Rot} \\ \M{R}_{H+Rot} \end{array}\right]^{\sf{T}}
\label{eq2027}
\end{equation}
where $\M{R}_{H+Rot}$ is the enriched coarse space and 
\begin{equation}\label{eq:juan002c}
\M{K}_{E,0}=\left[\begin{array}{c} \M{R}_{H+Rot} \\ \M{R}_{H+Rot} \end{array}\right] \M{K}_E \left[\begin{array}{c} \M{R}_{H+Rot} \\ \M{R}_{H+Rot} \end{array}\right]^{\sf{T}}.
 \end{equation}

\subsection{Elasticity operator local problems and coarse basis from a diffusion operator}\label{sec:53}
The diagonal blocks associated with every coarse partition are relatively small. Even though they are three times larger compared to the diagonal blocks obtained for the Laplace problem, the computational time does not increase significantly. Thus the diagonal blocks of the first part of the preconditioner can be obtained directly from the elastic operator. The second part can be constructed similar to the previous two cases using a coarse space from the diffusion operator and enriched coarse space with rotations. That is we have, 
\begin{equation}\label{eq:juan001}
\M{M}_{EH}^{-1}= \M{M}^{-1}_{E,1}+ 
\left[\begin{array}{c} \M{R}_{H}  \\ \M{R}_H \end{array}\right] \M{K}_{E,0}^{-1}  \left[\begin{array}{c} \M{R}_{H} \\ \M{R}_H \end{array}\right]^{\sf{T}}
\end{equation}
and 
\begin{equation}\label{eq:juan002}
\M{M}_{EH+Rot}^{-1}= \M{M}^{-1}_{E,1}+ 
\left[\begin{array}{c} \M{R}_{H+Rot}  \\ \M{R}_{H+Rot} \end{array}\right] \M{K}_{E,0}^{-1}  \left[\begin{array}{c} \M{R}_{H+Rot} \\ \M{R}_{H+Rot} \end{array}\right]^{\sf{T}}
\end{equation}
where we use the definitions in \eqref{eq:juan001c}
and  \eqref{eq:juan002c}, respectively. 
The first level of the preconditioner,  $\M{M}^{-1}_{E,1}$,  
is defined in \eqref{eqjuan:localproblemsE}.
These two cases complete the definitions of the preconditioners utilized in the numerical experiments. See Table \ref{tab:methods} for a summary of all the implemented iterations.
\begin{table}
\begin{tabular}{l l l l l}
\toprule
Notation &{Definition}  & {Level 1}  & {Eigenproblem}& Enrichment\\\midrule\midrule
$M_{EE}$  &  Eq. \eqref{eq:2lEp}  &  Elasticity  & Full elasticity \eqref{eq:juanelasticityeigenproblem} & None \\\midrule
$M_{HH}$  &  Eq. \eqref{eq2025} & Heat (blocks)   & Full heat \eqref{eq:juaneheateigenproblem} & None\\
$M_{HH+Rot}$  &  Eq. \eqref{eq2027} & Heat (blocks)  & Full heat \eqref{eq:juaneheateigenproblem} & Rotations \\\midrule
 $M_{EH}$ &  Eq. \eqref{eq:juan001} & Elasticity & Full heat \eqref{eq:juaneheateigenproblem} & None \\
$M_{EH+Rot}$  &  Eq. \eqref{eq:juan002} & Elasticity & Full heat \eqref{eq:juaneheateigenproblem} & Rotations \\\midrule
$M_{EH+Rot;Rand}$  &  Similar to Eq. \eqref{eq:juan002} & Elasticity  & Randomized heat \eqref{eq:juaneheateigenproblem} & Rotations\\
$M_{EE;Rand}$  &  Similar to Eq. \eqref{eq:2lEp} & Elasticity  & Randomized elasticity \eqref{eq:juanelasticityeigenproblem}& None
\\
\bottomrule
\end{tabular}
\caption{Information about the elasticity preconditioners used for numerical tests. 
All the methods solve an elasticity coarse problem, the differences can be found in the local solvers (Column 3) and in the construction of the coarse space where the changes correspond to the neighborhood eigenvalue problem used and its numerical approximation (Column 4). Some methods need enrichment of the coarse space (Column 5).
}
\label{tab:methods}
\end{table}

\section{Numerical experiments and results}
\label{sec::num}
\begin{figure}
\centering
\includegraphics[width=.3\textwidth]{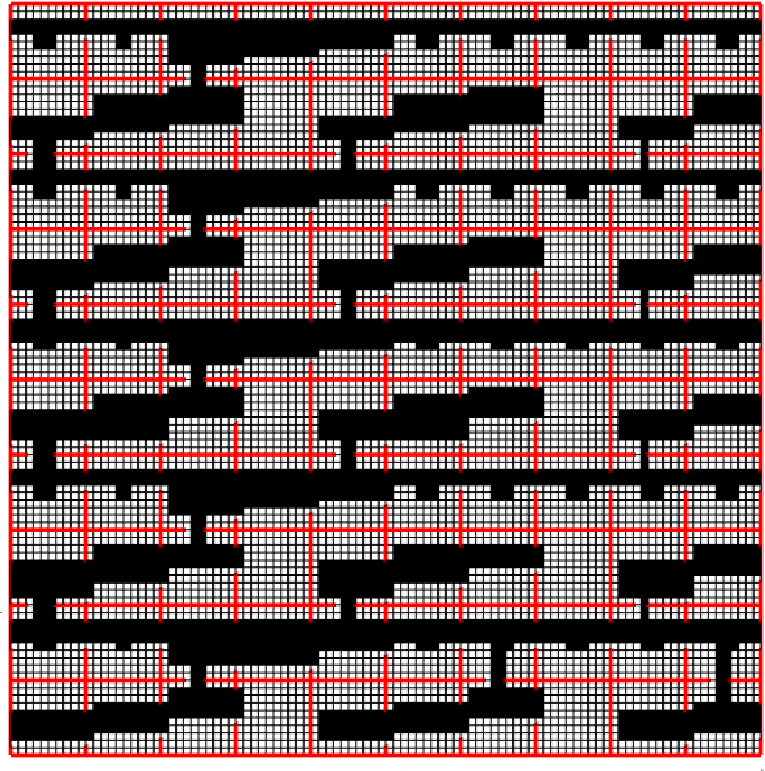}
\caption{Material distribution and coarse mesh for the numerical experiments.}
\label{fig:highcontrast}
\end{figure}

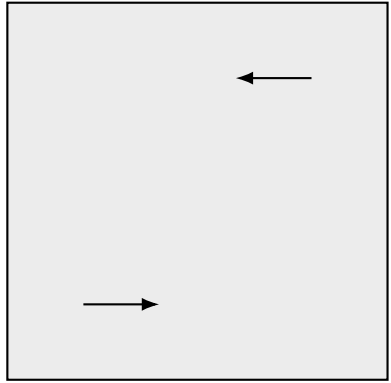
\begin{figure} 
    \centering
    \begin{tikzpicture}[scale=5]
        \fill [fill=gray!15!white] (0,0) rectangle (1,1);
        \draw [black,thick] (0,0) rectangle (1,1);
        \draw [-Latex, thick] (0.2,0.2) to (0.2+0.2,0.2);
        \draw [-Latex, thick] (0.8,0.8) to (0.8-0.2,0.8);
    \end{tikzpicture}
    \caption{Forcing term for the numerical experiments.}
    \label{fig:loading}
\end{figure}

In addition to the topology optimized case discussed earlier, in order to demonstrate the contrast independence of the above preconditioners, we utilize high-contrast material distribution as shown in \autoref{fig:highcontrast}. The computational domain is partitioned applying $10\times 10$ coarse mesh and each coarse-element is further partitioned utilizing  $10\times 10$ fine-mesh. The fine discretization is performed with bilinear polynomials. Zero Dirichlet boundary conditions are applied to all boundaries of the computational domain. The tests are performed with forcing terms shown in \autoref{fig:loading}. Both forces are applied to the solid-material subdomain which is the expected case for topology optimization problems. 
For the solution of the linear system, we run PCG until the relative norm of the initial residual is reduced by a factor of $10^{-6}$. 
We use the Lanczos connection method to estimate the condition number of the preconditioned operator; see \cite[Chapter 6]{saad2003iterative}.

In Table \ref{tab:methods} we present the results with the implemented preconditioners for the elasticity equation. 
After the construction of the coarse basis functions, the second level of all the implemented methods consists in  solving an elasticity coarse problem.  Some of the preconditioners differ in the
level-one local solvers: they either use elasticity equation local solvers or block diagonal heat equation local solvers. 
We then have several proposed coarse spaces constructions. 
The different  constructions correspond to  the local eigenvalue problem used and its numerical approximation (Column 4). We pose either a local elasticity eigenvalue problem or a local diffusion eigenvalue problem. For the numerical approximation of the eigenvalues and eigenvectors 
we use either a full eigenvalue solver of the fine scale local eigenvalue problem or the randomized method of Subsection \ref{rndalg}.
Some coarse spaces need to be  enriched (in order to obtain the RBM) by adding rotations  multiplied by partition of unity functions as additional coarse basis functions (Column 5).

In Table \ref{tab:elasticidad1}-\ref{tab:elasticidad6} we present iteration count and condition number estimates for the different iterations that have been introduced as they behave with respect to the contrast. We summarized this results in Table \ref{tab:num.exp.iterations} and \ref{tab:num.exp.spectcond}. We observe that the methods that are robust with respect to the contrast are $M_{EH+Rot}$, $M_{EE+Rot}$, 
$M_{EH+Rot; Rand}$ and $M_{EEt; Rand}$; see Table \ref{tab:methods}. 
Therefore,  they are computationally efficient alternatives to solve the elasticity equation and can be used in topology optimization problems such as the one considered in this paper. As it was mentioned before, an advantage of $M_{EH+Rot}$ is that, according to our numerical experiments and for general multiscale configurations, it allows us to identify the contrast-asymptotically-vanishing eigenvalues and the corresponding eigenvectors.

\begin{table}
\centering
\begin{tabular}{P{6cm} r r r r}
\toprule
{Preconditioner} & {Iterations} & {Condition} & Coarse dim.\\\midrule\midrule

{None} & 293 & $3.2\times 10^3$ & {---} \\\midrule 
$M_{EE}$& 14 & 4.6 & 243 \\
$M_{EE}$ using twice many eigenvectors& 15 & 5.0 & 486 \\\midrule
$M_{HH}$ & 29 & 21.0 & 243 \\
$M_{HH+Rot}$ & 29 & 21.0 & 243 \\\midrule
$M_{EH}$ & 15 & 5.2 & 243 \\
$M_{EH+Rot}$& 14 & 4.6 & 243 \\\midrule
$M_{EH+Rot;Rand}$ with 10 snapshots & 24 & 13.8 & 243 \\
$M_{EH+Rot;Rand}$ with 15 snapshots  & 25 & 15.6 & 243 \\
$M_{EE;Rand}$ with 10 snapshots & 20 & 9.4 & 243 \\
$M_{EE;Rand}$ with 15 snapshots & 20 & 9.4 & 243 \\
\bottomrule
\end{tabular}
\caption{Results for the elasticity problem with contrast $\eta=10^0$. See Table \ref{tab:methods} for the description of 
the preconditioners. }
\label{tab:elasticidad1}
\end{table}

\begin{table}
\centering
\begin{tabular}{P{6cm} r r r r}
\toprule
{Preconditioner} & {Iterations} & {Condition} & {Coarse dim.}\\\midrule

{None} & 1583 & $1.2\times 10^5$ & {---} \\\midrule
$M_{EE}$ & 26 & 16.0 & 387 \\
$M_{EE}$ using twice many eigenvectors & 19 & 7.7 & 774 \\\midrule
$M_{HH}$   & 86 & $2.4\times 10^2$ & 387 \\
$M_{HH+Rot}$   & 81 & $42.4\times 10^2$ & 387 \\\midrule
$M_{EH}$ & 35 & 72.0 & 387 \\
$M_{EH+Rot}$ & 28 & 18.3 & 387 \\\midrule
$M_{EH+Rot;Rand}$ with 10 snapshots & 32 & 26.8 & 387 \\
$M_{EH+Rot;Rand}$ with 15 snapshots & 33 & 27.0 & 387 \\
$M_{EE;Rand}$ with 10 snapshots & 30 & 22.9 & 387 \\
$M_{EE;Rand}$ with 15 snapshots & 30 & 22.8 & 387 \\
\bottomrule
\end{tabular}
\caption{Results for the elasticity problem with contrast $\eta=10^2$. See Table \ref{tab:methods} for the description of 
the preconditioners.}
\label{tab:elasticidad2}
\end{table}

\begin{table}
\centering
\begin{tabular}{P{6cm} r r r r}
\toprule
{Preconditioner} & {Iterations} & {Condition} & {Coarse dim.}\\\midrule
None & {\textgreater 2000} & $1.2\times 10^6$ & {---} \\\midrule
$M_{EE}$ & 53 & 113.8 & 387 \\
$M_{EE}$ using twice many eigenvectors & 28 & 23.6 & 774 \\\midrule
$M_{HH}$  & 200 & $2.1\times 10^3$ & 387 \\
$M_{HH+Rot}$  & 117 & $6.3\times 10^2$ & 387 \\\midrule
$M_{EH}$& 69 & $3.6\times 10^2$ & 387 \\
$M_{EH+Rot}$ & 56 & 108.9 & 387 \\\midrule
$M_{EH+Rot;Rand}$ with 10 snapshots  & 62 & 111.1 & 387 \\
$M_{EH+Rot;Rand}$ with 15 snapshots  & 62 & 111.4 & 387 \\
$M_{EE;Rand}$ with 10 snapshots& 57 & 121.9 & 387 \\
$M_{EE;Rand}$ with 15 snapshots & 58 & 122.0 & 387 \\
\bottomrule
\end{tabular}
\caption{Results for the elasticity problem with contrast $\eta=10^4$. See Table \ref{tab:methods} for the description of 
the preconditioners.}
\label{tab:elasticidad4}
\end{table}

\begin{table}
\centering
\begin{tabular}{P{6cm} r r r r}
\toprule
{Preconditioner} & {Iterations} & {Condition} & {Coarse dim.}\\\midrule
{None} & {\textgreater 2000} & $4\times 10^6$ & {---} \\\midrule
$M_{EE}$ & 44 & 140.6 & 387 \\
$M_{EE}$ using twice many eigenvectors & 26 & 15.1 & 774 \\\midrule
$M_{HH}$ & 141 & $2.2\times 10^4$ & 387 \\
$M_{HH+Rot}$ & 110 & $5.4\times 10^2$ & 387 \\\midrule
$M_{EH}$& 58 & 1.8e2 & 387 \\
$M_{EH+Rot}$& 58 & 140.8 & 387 \\\midrule
 $M_{EH+Rot;Rand}$ with 10 snapshots & 69 & 276.7 & 387 \\
 $M_{EH+Rot;Rand}$ with 15 snapshots & 69 & 276.5 & 387 \\
$M_{EE;Rand}$ with 10 snapshots & 63 & 141.0 & 387 \\
$M_{EE;Rand}$ with 15 snapshots  & 63 & 141.1 & 387 \\
\bottomrule
\end{tabular}
\caption{Results for the elasticity problem with contrast $\eta=10^6$. See Table \ref{tab:methods} for the description of 
the preconditioners.}
\label{tab:elasticidad6}
\end{table}

\begin{table}
\centering
\begin{tabular}{lrrrr}
&\multicolumn{4}{c}{\textbf{Contrast}}\\
\toprule
\multicolumn{1}{c}{\textbf{Preconditioner}} & \multicolumn{1}{c}{1} & \multicolumn{1}{c}{$1\times 10^{-2}$}& \multicolumn{1}{c}{$1\times 10^{-4}$} & \multicolumn{1}{c}{$1\times 10^{-6}$} \\\midrule
{None} & 292 & 1583 & {$>2000$} & {$>2000$} \\
$M_{EE}$ & 14 & 26 & 53&44 \\
$M_{EH+Rot}$ & 14 & 28 & 56 & 58 \\
$M_{EH+Rot;Rand}$ & 25 & 33 & 62 & 69 \\
$M_{EE;Rand}$ & 20 & 30 & 58 & 63 \\
\bottomrule
\end{tabular}
\caption{PCG iterations for different contrast values. See Table 
\ref{tab:methods} for description of the preconditioners.}
\label{tab:num.exp.iterations}
\end{table}

\begin{table}
\centering
\begin{tabular}{lrrrr}
&\multicolumn{4}{c}{\textbf{Contrast}}\\
\toprule
\multicolumn{1}{c}{\textbf{Preconditioner}} & \multicolumn{1}{c}{1} & \multicolumn{1}{c}{$1\times 10^{-2}$}& \multicolumn{1}{c}{$1\times 10^{-4}$} & \multicolumn{1}{c}{$1\times 10^{-6}$} \\\midrule
None & 3.2e3 & 1.2$\times 10^6$ & 1.2$\times 10^6$ & 4$\times 10^6$ \\
$M_{EE}$ & 4.6 & 16 & 113.8 & 140.6 \\
$M_{EH+rot}$ & 4.6 & 18.3 & 140.8 & 140.8 \\
$M_{EH+rot; Rand}$ & 15.6 & 27 & 111.4 & 276.5\\
$M_{EE; Rand}$ & 9.4 & 22.8 & 122 & 141.1 \\
\bottomrule
\end{tabular}
\caption{Spectral condition number of $M^{-1}A$ for different contrast values.  See Table \ref{tab:methods} for the description of 
the preconditioners.}
\label{tab:num.exp.spectcond}
\end{table}%

\newpage
\section{Conclusions} \label{S:concl}
In this paper we design and implement robust (with respect to the high-contrast and the multiscale structure) two-levels domain decomposition preconditioners for the elasticity equation appearing in topology optimization problems. Our design fits within the framework of the GMsFEM methodology where approximations of locally posed eigenvalue problems are used to construct the coarse space. We present several low cost  constructions with similar number of iterations. The computational cost related to the construction of the preconditioners is reduced an order of magnitude. The presented numerical experiments demonstrate the quality and robustness of our iterations. 




\bibliography{bib}

\end{document}